\catcode`\@=11
\magnification=1200 \parindent=0mm  \hfuzz=10pt
\hsize=13cm  \vsize=19cm   \hoffset=4mm		 \voffset=1cm
\pretolerance=500   \tolerance=1000   \brokenpenalty=5000
\catcode`\;=\active
\def;{\relax\ifhmode\ifdim\lastskip>\z@
\unskip\fi\kern.2em\fi\string;}
\catcode`\:=\active
\def:{\relax\ifhmode\ifdim\lastskip>\z@\unskip\fi
\penalty\@M\ \fi\string:}
\catcode`\!=\active
\def!{\relax\ifhmode\ifdim\lastskip>\z@
\unskip\fi\kern.2em\fi\string!}
\catcode`\?=\active
\def?{\relax\ifhmode\ifdim\lastskip>\z@
\unskip\fi\kern.2em\fi\string?}
\def\^#1{\if#1i{\accent"5E\i}\else{\accent"5E #1}\fi}
\def\"#1{\if#1i{\accent"7F\i}\else{\accent"7F #1}\fi}
\frenchspacing
\def\p{\par}
\def\x{\times}
\def\hs{\hskip}
\def\vs{\vskip}
\def\fin{$\diamondsuit$}
\catcode`\@=11
\def\system#1{\left\{\null\,\vcenter{\openup1\jot\m@th
\ialign{\strut\hfil$##$&$##$\hfil&&\enspace$##$\enspace&
\hfil$##$&$##$\hfil\crcr#1\crcr}}\right.}
\catcode`\@=12
\def\Fhd#1#2{\smash{\mathop{\hbox to 14mm{\rightarrowfill}}
\limits^{\scriptstyle#1}_{\scriptstyle#2}}}
\def\Fhg#1#2{\smash{\mathop{\hbox to 14mm{\leftarrowfill}}
\limits^{\scriptstyle#1}_{\scriptstyle#2}}}
\def\fhd#1#2{\smash{\mathop{\hbox to 8mm{\rightarrowfill}}
\limits^{\scriptstyle#1}_{\scriptstyle#2}}}
\def\fhg#1#2{\smash{\mathop{\hbox to 8mm{\leftarrowfill}}
\limits^{\scriptstyle#1}_{\scriptstyle#2}}}

\def\fvb#1#2{\llap{$\scriptstyle#1$}\left\downarrow
\vbox to 4mm{}\right.\rlap{$\scriptstyle#2$}}
 
\def\diagram#1{\def\normalbaselines{\baselineskip=0pt
\lineskip=10pt\lineskiplimit=1pt}   \matrix{#1}}
\count10=0  \count11=1  \count12=1 \count13=1 \count14=1
\def\pr#1{\if#11{^{'}}\else{\if#12{^{"}}\else{^{(#1)}}\fi}\fi}
  \def\parno{\count11}  
\def\bibno{\count14}

\def\bi#1#2#3{\item{{\bf [\the\bibno]}}
{\bf #1}{\it #2}{#3}\vskip 3mm plus 1mm minus 1mm \advance
\bibno by 1}
\def\tprg#1{\advance\parno by 1\vskip 2mm plus1mm minus 1mm
\noindent {\bf \the\parno  #1}
\count12=1 \count13=1\nobreak $$ $$}

\centerline{\bf Equivalence de la th\'eorie homotopique des
$n$-groupo\"{i}des}
\centerline{\bf et celle des espaces topologiques $n$-tronqu\'es }
\p\vs   1cm
\centerline{\bf Zouhair TAMSAMANI}
\p\vs   1cm
{\bf Introduction :}
\p\vs 4mm\hs 4mm
La pr\'esente note est une continuation d'un periprint, sur les
$n$-cat\'egories et $n$-groupo\"{i}des
non stricte [Ti], en vue de montrer le Th\'eor\`eme suivant :
\p\vs 2mm
{\bf Th\'eor\`eme (3.0) :} {\it On d\'esigne par $n$-Ho-Top la
cat\'egorie localis\'ee
des espaces topologiques $n$-tronqu\'es par les \'equivalences
d'homotopies faibles et
par $n$-Ho-Gr la cat\'egorie localis\'ee des $n$-groupo\"{i}des
par les $n$-\'equivalences
ext\'erieures. Alors les deux foncteurs
\p
\centerline{$\diagram{n-Ho-Top\ \Fhd{\Pi_{n}(\ )}{\Fhg{}{\mid	
\ \mid}}\ n-Ho-Gr}$}
\p
r\'ealisation g\'eom\'etrique $\mid\ \mid$ et $n$-groupo\"{i}de de
Poincar\'e $\Pi_{n}(\ )$ induisent une \'equivalence entre les deux
cat\'egories $n$-Ho-Top et $n$-Ho-Gr.}
\p\vs 2mm\hs 4mm
On aborde le probl\`eme en deux \'etapes :
la premi\`ere consiste \`a construire, pour chaque espace topologique
$n$-tronqu\'e $X$, un isomorphisme $\mid\Pi_{n}(X)\mid\fhd{}{}X$
dans la cat\'egorie $n$-Ho-{\cal T}op compatible avec la fonctorialit\'e
en $X$.
La deuxi\`eme est de construire, pour tout $n$-groupo\"{i}de $\Phi$, une
$n$-\'equivalence
ext\'erieure $\Phi\fhd{}{}\Pi_{n}(\mid\Phi\mid)$ compatible aussi avec
la fonctorialit\'e
en $\Phi $. Pour cela on s'inspire d'un
r\'esultat de Segal [Sl] sur les espaces topologiques simpliciaux qui
r\'esout d\'ej\`a
la deuxi\`eme \'etape pour $n=1$. On utilise aussi quelques r\'esultats de
Bousfield-Kan [BK] qui traitent aussi le cas $n=1$ et montre que la limite
homotopique directe d'un espace topologique bisimplicial est faiblement
\'equivalente
\`a la diagonale de celui-ci. Puis on introduit deux foncteurs $\Omega $ et
${\cal R}$ comme suite :
\p
\centerline{$\diagram{(n+1)-ETS\Fhd{\Omega }{\Fhg{}{{\cal R}}}n-ETS}$}
\p
qui permettent de donner une autre interpr\'etation au $n$-groupo\"{i}de de
Poincar\'e
et la r\'ealisation g\'eom\'etrique d'un espace topologique $n$-simplicial.
On montre, pour tout espace topologique $X$, on a le
diagramme d'\'equivalences d'homotopie faibles suivant :
\p
\centerline{$\diagram{{\cal R}^{n}\Omega^{n}X&\Fhg{\sim}{}&
{\cal R}^{n}(F\circ\Omega^{n}X)&\Fhd{\sim}{}
&{\cal R}^{n}(\pi_{0}\circ\Omega^{n}X)
\cr \fvb{\wr}{}&&&&\fvb{}{\wr}\cr X&&&&\mid\Pi_{n}(X)\mid}$}
\p
o\`u  $F(Y)=Y\x _{\pi_{0}(Y)^{g}}\pi_{0}(Y)^{d}$  paragraphe (3.3).
Lorsqu'on regarde ce diagramme dans la cat\'egorie $n$-Ho-{\cal T}op
on obtient un isomorphisme de $\mid\Pi_{n}(X)\mid$ vers $X$
compatible avec la fonctorialit\'e en $X$.
Au deuxi\`eme \'etape on construit pour tout $n$-groupo\"{i}de $\Phi $ une
$n$-\'equivalence ext\'erieure de $\Phi $ vers $\Pi_{n}(\mid \Phi \mid)$
compatible aussi avec la fonctorialit\'e en $\Phi $. D'autre part
on v\'erifie
que les foncteurs $\mid\ \mid$ et $\Pi_{n}(\ )$ pr\'eservent les
\'equivalences
par lesquelles on localise, ce montre que ces foncteurs induisent
effectivement
une \'equivalence entre les cat\'egories $n$-Ho-Top et $n$-Ho-Gr.
\p\hs 2mm
On appelle {\it espace topologique $n$-simplicial} un foncteur
contravariant
de la cat\'egorie $\Delta ^{n}$ vers celle des espaces topologiques.
Une transformation (resp. isomorphisme) naturelle entre de tels
foncteurs
sera appel\'ee {\it morphisme (resp. isomorphisme)}. Lorsque $n=0$
on retrouve
les espace topologiques et les applications continues et pour $n=1$ on dit
simplement espace topologique simplicial.
Soient $n\geq 0$ un entier et $X$ un espace topologique. On dit que $X$ est
{\it $n$-tronqu\'e} si et seulement si pour tout $x\in X$ les groupes
d'homotopies
$\pi_{i}(X,x)$ sont nuls pour $i > n$.
\p\vs 4mm
{\bf (3.1).--- QUELQUES R\'ESULTATS DE G. Segal ET Bousfield-Kan :}
\p\vs 4mm\hs 4mm
{\bf D\'efinition (3.1.1) :} {On appelle {\it espace simplicial de Segal}
un espace topologique simplicial $\Phi$ v\'erifiant les trois axiomes :
\p
(i) L'espace $\Phi_{0}$ est muni de la topologie discr\`ete.
\p
(ii) Pour tout entier $m \geq 2$, le morphisme :
\p\vs 4mm
\centerline{$\Phi_{m}\ \Fhd{\delta_{[m]}}{}\
\Phi_{1}{\x}_{\Phi_{0}}\cdots{\x}_{\Phi_{0}}\Phi_{1}$}}
\p\vs 2mm
est une \'equivalence d'homotopie faible (e.h.f).
\p
(iii) L'ensemble simplicial $\pi_{0}\circ\Phi$ est le nerf d'un
groupo\"{i}de.
(Les axiomes (i) et (ii) entra\^{i}nent que $\pi_{0}\circ\Phi$ le nerf
d'une cat\'egorie.)\p\vs 4mm\hs 4mm
La construction de la r\'ealisation g\'eom\'etrique d'une
$n$-pr\'e-cat\'egorie dans [11]
garde aussi un sens pour les espaces topologiques $n$-simpliciaux,
 mais dans
ce cas on tiendra compte de la topologie des espaces $\Phi_{M}$ o\`u
$M$ est un objet de
$\Delta ^{n}$. On rappelle que cette r\'ealisation g\'eom\'etrique
est d\'efinie par :\p
$$\mid\Phi \mid =[\displaystyle\coprod_{M}\Phi_{M}\x R^{M}]^{\sim}$$
\p
o\`u $\sim$ est la relation qui permet d'identifier l'\'el\'ement
$(x,\sigma ^{'}(a))$ \`a l'\'el\'ement $(\sigma ^{"}(x),a)$,
pour tout
$\sigma : M\fhd{}{}M^{'}$ et tout $(x,a)$ dans $\Phi_{_{M^{'}}}\x R^{^{M}}$.
La topologie de $\mid\Phi \mid$ sera induite de la topologie produit puis la
topologie
quotient. Lorsque $\Phi_{M}$ n'est pas munit d'une topologie, on lui
donne la topologie discr\`ete. \p\hs 4mm
Soient $\Phi$ un espace simplicial et $x$, $y$ deux \'el\'ements de
$\Phi_{0}$,
on d\'esigne par $\Phi_{1}(x,y)$ l'image inverse de $(x,y)$ par
le morphisme :
\p\vs 4mm
\centerline{$\Phi_{1}\ \Fhd{\delta^{'}_{0}\x\delta^{'}_{1}}
{}\ \Phi_{0}\x\Phi_{0}$}
\p\vs 4mm
\centerline{$\Phi_{1}(x,y)=\{f\in\Phi_{1}\ \mid\ \delta^{'}_{0}(f)=x\
et\ \delta^{'}_{1}(f)=y\}$}\p\vs 4mm
D'autre part on a une injection $\Phi_{0}\fhd{}{}\mid\Phi \mid$ qui
envoie tout \'el\'ement
$x$ de $\Phi_{0}$ vers la classe $\mid (x,1)\mid$ dans $\mid\Phi\mid$.
Dans la suite on
identifie $\Phi_{0}$ avec son image dans $\mid\Phi \mid$.
On a un morphisme : \p\vs 4mm
\centerline{$\Phi_{1}(x,y)\ \Fhd{\alpha}{}\ Ch^{x,y}({\cal R}\Phi )$}
\p\vs 4mm
La notation $Ch^{x,y}(\mid\Phi\mid )$ d\'esigne l'espace des chemins de
$x$ vers $y$
dans $\mid\Phi \mid$ et le morphisme $\alpha$ est ce lui qui \`a chaque
$f\in \Phi_{1}(x,y)$
fait correspondre le chemin $\Gamma : R\fhd{}{}\mid\Phi \mid$ tel que
$\Gamma (t)=\mid(f,t)\mid$. Ce chemin va bien de $x$ vers $y$ en effet :
\p
soient les morphismes $\delta^{"}_{0}, \delta^{"}_{1} : \{1\}\fhd{}{}R$ et
$\delta^{'}_{0}, \delta^{'}_{1} : \Phi_{1}\fhd{}{}\Phi_{0}$ alors on a :
\p\vs 4mm
\centerline{$(f,\delta^{"}_{0}(1))\sim(\delta^{'}_{0}(f),1)=(x,1)=x$}
\p\vs 4mm
\centerline{$(f,\delta^{"}_{1}(1))\sim(\delta^{'}_{1}(f),1)=(y,1)=x$}
\p\vs 4mm
Lorsque $\Phi $ est un espace topologique simplicial on notera
${\cal R}\Phi =\mid\Phi \mid$. On peut maintenant \'enoncer le th\'eor\`eme
principal de G. Segal.
\p\vs 4mm\hs 4mm
{\bf Th\'eor\`eme (3.1.2)} {\it Si $\Phi $ est un espace simplicial
de Segal, alors on a  :
\p\vs 2mm
(i) Le morphisme $\alpha : \Phi_{1}(x,y)\fhd{}{}Ch^{x,y}({\cal R}\Phi )$
est une e.h.f.
\p
(ii) L'injection $\Phi_{0}\fhd{}{}{\cal R}\Phi $ induit une bijection
du tronqu\'e $T(\pi_{0}\circ\Phi)$ sur l'ensemble
$\pi_{0}({\cal R}\Phi )$.}
\p\vs 4mm\hs 4mm
{\bf Preuve :} G. Segal [Sl]
\p\vs 4mm\hs 4mm
Faisons maintenant une construction dans une autre direction. Soit $X$
un espace topologique et d\'esignons par $Hom^{-}(R^{m},X)$ le produit
fibr\'e suivant :
\p\vs 2mm
\centerline{$Hom^{-}(R^{m},X):=Hom(R^{m},X)\x_{X^{m+1}}(X^{disc})^{m+1}$}
\p\vs 2mm
La notation $()^{disc}$ signifie le m\^{e}me espace avec la topologie
discr\`ete.
L'espace $Hom(R^{m},X)$ est munit de la topologie compact-ouvert.
L'application de $Hom(R^{m},X)$ vers $X^{m+1}$ est celle qui associe
\`a chaque
$m$-simplexe $f$ de $Hom(R^{m},X)$ ces ($m$+1) sommets :
\p\vs 2mm
\centerline{$(f(\delta^{"}_{i}))_{0\leq i\leq m}$ \hs 1cm o\`u \hs 1cm
$\delta^{"}_{i} : \{1\}\fhd{}{}R^{m}$}
\p\vs 2mm\hs 4mm
Soit $X$ un espace toplogique, on note par $\Omega X$ l'ensemble simplicial
d\'efini par :
\p\vs 2mm
\centerline{$(\Omega X)_{_{m}}:=Hom^{-}(R^{m},X)$}
\p\vs 2mm
Remarquons ici que pour $x,y$ dans $X$ on a :
$(\Omega X)_{_{1}}(x,y)=Ch^{x,y}(X)$.
\p\vs 4mm\hs 4mm
{\bf Proposition (3.1.3)} {\it Si $X$ est un espace topologique alors
$\Omega X$
est un espace simplicial de Segal.}
\p\vs 4mm\hs 4mm
{\bf Preuve :} Remarquons d'abord que $(\Omega X)$ est un espace
topologique
simplicial, tel que $(\Omega X)_{_{0}}=X^{disc}$. Consid\'erons les
morphismes
sommets et c\^{o}t\'es suivants :
\p\vs 4mm
\centerline{$\{1\}\ \Fhd{\delta^{"}_{i}}{}\ R^{m}\ \Fhg{\delta^{"}_{ij}}
{}\ R$\hs
1cm et posons \hs 1cm
$R_{i,i+1}=\delta^{"}_{i,i+1}(R)$ et $\delta^{"}_{i}(1)=\varepsilon^{i}$}
\p\vs 2mm
On sait qu'il existe une r\'etraction par d\'eformation $r$ de $R^{m}$ sur
la r\'eunion
de ces $m$ sommets :
\p\vs 3mm
\centerline{$R^{m}\Fhd{r}{}W:=
\displaystyle\bigcup _{0\leq i\leq m-1}R_{i,i+1}\Fhd{j}{}R^{m}$}
\p\vs 3mm
tels que $j\circ r \sim id_{R^{m}}$ et $r\circ j=id_{W}$. Le r\'etracte
$r$ induit
un inverse homotopique du morphisme :
\p\vs 2mm
\centerline{$(\Omega X)_{_{m}}\Fhd{\delta_{[m]}}{}(\Omega X)_{_{1}}
{\x}_{(\Omega X)_{_{0}}}
\cdots{\x}_{(\Omega X)_{_{0}}}(\Omega X)_{_{1}}$\hs 4mm par le morphisme}
\p\vs 4mm
\centerline{$\matrix{(\Omega X)_{_{1}}{\x}_{(\Omega X)_{_{0}}}
\cdots{\x}_{(\Omega X)_{_{0}}}(\Omega X)_{_{1}}\Fhd{\phi}
{}(\Omega X)_{_{m}}\cr\cr
\Bigl(\lambda_{i,i+1},(a_{i},a_{i+1})\Bigr)_{i}\Fhd{}{}
\Bigl(\lambda\circ r,(a_{i})_{i}\Bigr)}$}
\p\vs 4mm
o\`u $\lambda$ est le recollement des $(\lambda_{i,i+1})$
suivant le diagramme :
\p
\centerline{$\diagram{R_{i,i+1}&\Fhd{}{}&W&
\Fhg{r}{}&R^{m}\cr \fvb{\lambda_{i,i+1}}{}&&\fvb{\lambda}
{}&&\fvb{}{\lambda\circ r}
\cr X&\Fhd{}{id_{X}}&X&\Fhd{}{id_{X}}&X}$}
\p
On v\'erifiera facilement qu'on a : \hs 5mm $\delta_{[m]}
\circ\phi=id_{\cal B}$
\hs 3mm et \hs 3mm $\phi\circ\delta_{[m]}\sim id_{\cal A}$ \hs 5mm o\`u
\p\vs 2mm
\centerline{${\cal A}=(\Omega X)_{_{m}}$\hs 5mm et\hs 5mm
${\cal B}=(\Omega X)_{_{1}}{\x}_{(\Omega X)_{_{0}}}
\cdots{\x}_{(\Omega X)_{_{0}}}(\Omega X)_{_{1}}$}
\p\vs 2mm
Si ${\cal H}$ est l'homotopie $j\circ r \sim id_{R^{m}}$, alors l'homotopie
$\phi\circ\delta_{[m]}\sim id_{\cal A}$ est d\'efinie explicitement
omme suit :
\p
\centerline{$\diagram{I\x{\cal A}&\Fhd{{\cal K}}{}&{\cal A}\cr
\Bigl(t,(f,(a_{i}))\Bigr)&\Fhd{}{}&
\Bigl(f\circ{\cal H}(t,\ ),f\circ{\cal H}(0,\varepsilon_{i})\Bigr)}$}
\hfill\fin\p
D'autre part, si $X$ est un espace toplogique, on a un morphisme :
\p
\centerline{$\diagram{{\cal R}(\Omega X)&\Fhd{}{}&X\cr
\mid((f,a_{i}),x)\mid&\Fhd{}{}&f(x)}$}
\p\hs 4mm
{\bf Th\'eor\`eme (3.1.4)} {\it Pour tout espace topologique $X$
le morphisme :
\p\vs 2mm
\centerline{${\cal R}(\Omega X)\ \Fhd{}{}\ X$}
\p\vs 2mm
est une \'equivalence d'homotopie faible.}
\p\vs 4mm\hs 4mm
{\bf Preuve :} (i) V\'erifions d'abord qu'on a un isomorphisme entre
les deux nerfs
$\pi_{0}(\Omega X)$ et $\Pi_{1}(X)$. Soient $m$ un objet de $\Delta $, et
$\Gamma=\Gamma_{1}\x\Gamma_{2}$ un chemin dans $(\Omega X)_{_{m}}$
entre deux \'el\'ements $(f,(x_{i}))$ et $(g,(y_{i}))$ de $(\Omega X)_{_{m}}$.
Alors
$\Gamma_{1}$ est une homotopie entre $f$ et $g$ dans $Hom(R_{m},X)$,
et $\Gamma_{2}$
un morphisme de l'intervalle $I$ vers $(X^{disc})^{m+1}$ qui est
n\'ecessairement
constant, car $I=[0,1]$ est connexe. Donc $\Gamma_{1}$ est une homotopie qui
fixe les sommets. On en d\'eduit que le morphisme premi\`ere projection
$Pr_{1} : (\Omega X)_{m}\fhd{}{}X_{m}$ induit un isomorphisme
\p\vs 4mm
\centerline{$\pi_{0}(\Omega X)_{_{m}}\ \Fhd{\eta^{1}_{X}(m)}
{}\ \overline{(X_{m})}$}
\p\vs 4mm
qui est en plus fonctorielle en $m$. D'o\`u
$\eta^{1}_{X} : \pi_{0}(\Omega X)\fhd{}{}\Pi_{1}(X)$ est un
isomorphisme naturel.
\p
(ii) L'ensemble $\pi_{0}(X)$ est \'egale \`a $T[\Pi_{1}(X)]$
lequel en bijection avec
$T[\pi_{0}(\Omega X)]$ en vertu de (i), donc d'apr\`es le Th\'eor\`eme
de Segal le morphisme
${\cal R}(\Omega X)\fhd{}{}X$ induit un isomorphisme sur les $\pi_{0}$.
\p
(iii) Soient $x,y$ deux points de $X$ et consid\'erons les deux morphismes :
\p
\centerline{$\diagram{(\Omega X)_{_{1}}(x,y)&
\Fhd{\alpha }{}&Ch^{x,y}({\cal R}(\Omega X))&\Fhd{\beta }{}&Ch^{x,y}(X)}$}
\p
o\`u  $\alpha $ est une e.h.f
d'apr\`es le Th\'eor\`eme de Segal, et $\beta\circ \alpha $ un isomorphisme
par identification. Donc $\beta $ est une e.h.f  et induit
par la suite des isomorphismes sur les groupes d'homotopie de
${\cal R}(\Omega X)$ et $X$ en degr\'e $\geq 1$.
\p\hfill\fin
\p\vs 4mm\hs 4mm
{\bf Th\'eor\`eme (3.1.5)} {\it Si $t : \Phi \fhd{}{}\Psi  $ est un
morphisme entre
espaces simpliciaux, tel que pour tout objet $m$ de $\Delta $ le morphisme
$t_{m} : \Phi_{m} \fhd{}{}\Psi  _{m}$ est une e.h.f.
Alors ${\cal R}t : {\cal R}\Phi \fhd{}{}{\cal R}\Psi  $ est une
\'equivalence
d'homotopie faible.}
\p\vs 4mm\hs 4mm
R\'esumons d'abord les r\'esultats de Bousfield-Kan [BK] qu'on
va utiliser pour d\'emontrer ce Th\'eor\`eme :
\p
(i) Si $X$ est un espace topologique, alors le morphisme naturel :
\p
$$\mid Sin(X)\mid\ \Fhd{\sim}{}\ X$$
\p
est une \'equivalence d'homotopie faible, o\`u $Sin(X)$ est le foncteur
singulier
d\'efini par $Sin(X)=Hom(R^{m},X)$ ([BK] p. 329).
\p
(ii) Si $f : \Phi\fhd{}{}\Psi$ est un morphisme d'ensembles 2-simpliciaux
tel que pour tout entier $m\geq 0$ le morphisme
$\mid f_{m}\mid :\mid \Phi_{m}\mid\fhd{}{}\mid\Psi_{m}\mid$ est une e.h.f,
 alors
$f$ induit une e.h.f sur les limites homotopiques
$\mid holim\Phi\mid\fhd{}{}\mid holim\Psi\mid$ ([BK] p. 329).
\p
(iii) Si $\Phi$ est un ensembles 2-simplicial, alors on a une e.h.f
$\mid holim\Phi\mid\fhd{}{}\mid diag\Phi\mid$ o\`u $diag\Phi$ est
la diagonale de
$\Phi$ d\'efinie par $(diag\Phi)_{m}=\Phi_{m,m}$ ([BK] p. 331).
\p
(iv) On peut v\'erifier facilement que pour tout ensemble simplicial $\Phi$
on a une e.h.f :
\p
$$\mid diag Sin\Phi\mid\fhd{}{}\mid\Phi\mid$$
\p\vs 4mm\hs 4mm
{\bf Preuve :} (1) Pour tout entier $m\geq 0$ et gr\^ace \`a (i) on a
un diagramme commutatif :
\p
\centerline{$\diagram{\mid Sin\Phi_{m}\mid&\Fhd{\mu}{}&\mid Sin\Psi_{m}\mid
\cr \fvb{\wr}{}&&\fvb{}{\wr}\cr
\Phi_{m}&\Fhd{\sim}{}& \Psi_{m}}$}
\p
qui induit un autre diagramme commutatif sur les $\pi_{i}$, donc $\mu$ est
une e.h.f.
\p
(2) En vertu de (ii), (iii) et (iv) on a les deux diagrammes commutatifs :
\p
\centerline{$\diagram{\mid holim Sin\Phi\mid&\Fhd{\sim}{}&\mid diag Sin\Phi\mid
&\Fhd{\sim}{}&\mid \Phi\mid
\cr \fvb{\wr}{}&&\fvb{}{\phi}&&\fvb{}{\mid f\mid}\cr
\mid holim Sin\Psi\mid&\Fhd{}{\sim}& \mid diag Sin\Psi\mid
&\Fhd{}{\sim}&\mid \Psi\mid}$}
\p
on en d\'eduit de m\^eme que $\phi$ est une e.h.f et puis $\mid f\mid$
est une e.h.f
\hfill\fin
\p\vs 4mm\hs 4mm
Maintenant on va proc\'eder \`a une construction de deux foncteurs qui
vont permettre
gr\^ace aux r\'esultats pr\'ec\'edents de donner une interpr\'etation
r\'ecurrente du
$n$-groupo\"{i}de de Poincar\'e $\Pi_{n}(X)$ d'un espace topologique
ainsi que sa
r\'ealisation g\'eom\'etrique construis dans [11]. En suite on va
utiliser
cette nouvelle construction pour montrer que les groupes d'homotopies
des espaces
$X$ et $\mid \Pi_{n}(X)\mid$ sont isomorphes.
\p\vs 4mm
{\bf (3.2).--- LES FONCTEURS ${\cal R}$ ET $\Omega $ :}
\p\vs 4mm\hs 4mm
Soit $n$ un entier positif. On d\'esigne par $n$-ETS la cat\'egorie
des espaces
topologiques $n$-simpliciaux et par $n$-ES celle des ensembles
$n$-simpliciaux.
\p
On d\'esigne par ${\cal R}$ : ($n$+1)-ETS $\fhd{}{}$ $n$-ETS le
foncteur covariant
d\'efinit pour tout $M$ objet de $\Delta ^{n+1}$ par
$({\cal R}\Phi )_{_{M}} :={\cal R}(\Phi _{_{M}})$. On obtient
alors une suite de
ce genre de foncteurs :
\p\vs 3mm
\centerline{$\fhd{{\cal R}}{}$ ($n$+1)-ETS $\fhd{{\cal R}}{}$ $n$-ETS
$\fhd{{\cal R}}{}$ \hs 1cm $\fhd{{\cal R}}{}$ ETS $\fhd{{\cal R}}{}$
{\cal T}op}
\p\vs 3mm
D'autre part on d\'efinit un autre foncteur covariant
$\Omega $ : $n$-ETS $\fhd{}{}$ ($n$+)-ETS pour tout $(M,m)$ objet de
$\Delta ^{n}\x\Delta $ par
$(\Omega \Phi )_{_{M,m}} := (\Omega\Phi _{_{M}})_{_{m}}$. De m\^eme
on a une suite de foncteurs :
\p\vs 3mm
\centerline{ {\cal T}op $\fhd{\Omega}{}$ ETS  $\fhd{\Omega }{}$ \hs 1cm
$\fhd{\Omega }{}$ $n$-ETS $\fhd{\Omega }{}$ ($n$+1)-ETS
$\fhd{\Omega }{}$}
\p\vs 4mm\hs 4mm
{\bf Proposition (3.2.1) :} {\it  Pour tout espace topologique $X$ il existe
un isomorphisme
naturel entre les deux $n$-pr\'e-cat\'egories $\pi_{0}(\Omega^{n}X)$ et
$\Pi_{n}(X)$
compatible avec la fonctorialit\'e en $X$. Et on en d\'eduit que
$\pi_{0}(\Omega^{n}X)$
est un $n$-groupo\"{i}de.}
\p\vs 4mm\hs 4mm
{\bf Preuve :} Montrons la proposition par r\'ecurrence sur $n$.
Soient $X$ un espace
topologique et $m$ un objet de $\Delta $. On sait d'apr\`es la
preuve du
Th\'eor\`eme (3.1.4) que le morphisme premi\`ere projection
$Pr_{1} : (\Omega X)_{m}\fhd{}{}X_{m}$ induit un isomorphisme naturel
$\eta^{1}_{X} : \pi_{0}(\Omega X)\fhd{}{}\Pi_{1}(X)$. On peut v\'erifier
facilement que
$\eta^{1}_{X}$ est compatible avec la fonctorialit\'e en $X$ des foncteurs
$\pi_{0}(\Omega X)$ et $\Pi_{1}(X)$.
La proposition est donc vraie pour $n=1$. Formulons notre hypoth\`ese
de r\'ecurrence comme suit : on suppose que jusqu'\`a un rang $n$ il
existe un isomorphisme
naturel $\eta^{n}_{X} : \pi_{0}(\Omega^{n}X)\fhd{}{}\Pi_{n}(X)$ compatible
avec la fonctorialit\'e
en $X$. Posons $Y=(\Omega X)_{m}$ alors il existe un isomorphisme
$\eta^{n}_{Y} : \pi_{0}(\Omega^{n}Y)\fhd{}{}\Pi_{n}(Y)$.
Comme $\Omega X$ est un foncteur, il s'en suit que $\eta^{n}_{Y}$ est
fonctoriel en $m$.
Pour $M$ objet de $\Delta ^{n}$, la projection $Pr_1$ induit un
isomorphisme :
\p\vs 2mm
$$ Y_{M}=(\Omega X)_{m,M}\ \Fhd{}{}\ X_{m,M}$$
\p\vs 2mm
 qui pr\'eserve les relations d'\'equivalences et induit \`a son tour
 un isomophpisme
\p\vs 2mm
$$F_{m,M} : \overline{ Y_{M}}\ \Fhd{}{}\ \overline{X_{m,M}}$$
\p\vs 2mm
fonctoriel en $X$ et en $(m,M)$. Par suite on obtient un isomorphisme
\p\vs 2mm
$$F_{m} : \Pi_{n}\Bigl((\Omega X)_{m}\Bigr)\ \Fhd{}{}\ \Pi_{n+1}(X)_{m}$$
\p\vs 2mm
Le compos\'e des deux isomorphismes naturels :
\p\vs 2mm
$$\pi_{0}(\Omega^{n}X)_{m}\ \Fhd{\eta^{n}_{Y}}{}\ \Pi_{n}
\Bigl((\Omega X)_{m}\Bigr)
\ \Fhd{F_{m}}{}\ \Pi_{n+1}(X)_{m}$$
\p\vs 2mm
est un isomorphisme fonctoriel en $X$ et en $m$. On en d\'eduit
un isomorphisme :
\p\vs 2mm
$$\pi_{0}(\Omega^{n}X)\ \Fhd{\eta^{n+1}_{X}}{}\ \Pi_{n+1}(X)$$
\p\vs 2mm
compatible avec la fonctorialit\'e en $X$.
\hfill\fin
\p\vs 4mm\hs 4mm
{\bf Remarque :} Un raisonnement analogue \`a celui utilis\'e dans la
preuve de la
proposition (3.2.1) montre qu'on a aussi, pour tout espace topologique
$X$ et tout
objet $M$ de $\Delta ^{n}$, un isomorphisme naturel :
\p\vs 2mm
$$\Pi_{1}\Bigl((\Omega^{n}X)_{M}\Bigr)\ \Fhd{}{}\ \Pi_{n+1}(X)_{M}$$
\p\vs 2mm
compatible avec la fonctorialit\'e en $X$ et en $M$.
\p\vs 4mm\hs 2mm
Nous avons construit pour chacun des $n$-espaces topologique simpliciaux
$\pi_{0}(\Omega^{n}X)$ et $\Pi_{n}(X)$ une r\'ealisation g\'eom\'etrique
de deux fa\c{c}on
diff\'erentes. Comme ces deux $n$-espaces sont isomorphes, il est naturel
de voir comment se comportent
leur r\'ealisations g\'eom\'etriques.
\p\vs 4mm\hs 4mm
{\bf Proposition (3.2.2) :} {\it Soit $\Phi $ un objet de $n$-ETS. Alors
les deux
espaces topologiques ${\cal R}^{n}\Phi $ et $\mid \Phi \mid$ sont
hom\'eomorphes.}
\p\vs 4mm\hs 4mm
{\bf Preuve :} Pour cela on va faire une r\'ecurrence sur $n$.
Lorsque $n=1$ on a
${\cal R}\Phi =\mid \Phi \mid$. Supposons que la proposition soit vraie
jusqu'\`a un rang $n$ et soit $\Phi $ un objet de ($n$+1)-ETS. En appliquant
l'hypoth\`ese de r\'ecurrence \`a ${\cal R}\Phi $, on obtient un
hom\'eomorphisme :
\p\vs 2mm
$${\cal R}^{n+1}\Phi \ \Fhd{}{}\ \mid {\cal R}\Phi \mid $$
\p\vs 2mm
On est donc amen\'e \`a montrer que $\mid {\cal R}\Phi \mid$ est
hom\'eomorphe
\`a $\mid \Phi \mid$. Pour tout $M$ objet de $\Delta$ on a un
morphisme naturel :
\p\vs 2mm
$$\displaystyle\coprod_{m}\Phi_{M,m}\x R^{m}\ \Fhd{}{}\
({\cal R}\Phi) _{M} $$
\p\vs 2mm
Et par suite un morphisme :
\p\vs 2mm
$$\matrix{{\cal A}=\displaystyle\coprod_{m,M}\Phi_{M,m}\x R^{m}
\x R^{M}\ \Fhd{\phi}{}\
\displaystyle\coprod_{M}({\cal R}\Phi )_{M}\x R^{M}={\cal B}\cr\cr
(a,(x_{0},x))\ \Fhd{}{}\ \Bigl(\mid (a,x_{0})\mid ,x\Bigr)}$$
\p\vs 2mm
ce morphisme pr\'eserve les relations d'\'equivalences d\'efinies sur
${\cal A}$ et ${\cal B}$, en effet :
\p
Soient $(a,(x_{0},x))$ un \'el\'ement de $\Phi_{M,m}\x R^{m^{'}}
\x R^{M^{'}}$ et
$\tau \x\sigma : [M]\x[m]\fhd{}{}[M^{'}]\x[m^{'}]$ une fl\`eche de
$\Delta^{n}\x\Delta$. Les deux \'el\'ements $(a,(\tau \x\sigma)^{"}(x_{0},x))$
et $((\tau \x\sigma)^{'}a,(x_{0},x))$ de ${\cal A}$ repr\'esentent
la m\^eme classe dans $\mid \Phi\mid $. Montrons alors que leur
images par $\phi$ repr\'esentent aussi la m\^eme classe dans
$\mid {\cal R}\Phi \mid$. On a :
\p\vs 2mm
$$\matrix{\phi \Bigl(a,(\tau \x\sigma)^{"}(x_{0},x)\Bigl)&=&
\Bigl(\mid(a,\sigma ^{"}x_{0})\mid,\tau ^{"}x)\Bigl)\cr\cr
\phi \Bigl((\tau \x\sigma)^{'}a,(x_{0},x)\Bigl)&=&
\Bigl(\mid((\tau \x\sigma)^{'}a,x_{0})\mid,x\Bigl)}$$
\p\vs 2mm
En utilisant les deux diagrammes commutatifs suivants :
\p
\centerline{$\diagram{[M]\x[m]&\fhd{I_{[M]}\x\sigma }{}&[M]\x[m^{'}]\cr
\fvb{\tau \x I_{[m]}}{}&&\fvb{}{\tau \x I_{[m^{'}]}}\cr [M^{'}]\x[m]&
\fhd{}{I_{[M^{'}]}\x\sigma }&[M^{'}]\x[m^{'}]}$\hs 2cm
$\diagram{\Phi _{_{M^{'},m^{'}}}&\fhd{(I_{[M^{'}]}\x\sigma )^{'}}{}&
\Phi _{_{M^{'},m}}\cr
\fvb{(\tau \x I_{[m^{'}]})^{'}}{}&&\fvb{}{(\tau \x I_{[m]})^{'}}\cr
\Phi _{_{M,m^{'}}}&\fhd{}{(I_{[M]}\x\sigma)^{'}}&\Phi _{_{M,m}}}$}
\p
on obtient les relations suivantes :
\p\vs 2mm
$$\matrix{\mid((\tau \x\sigma)^{'}a,x_{0})\mid&=&
\mid((\tau \x I_{[m]})^{'}\x (I_{[M^{'}]}\x\sigma )^{'}a,x_{0})\mid\cr\cr
&=&\mid\tau ^{'}\mid\mid((I_{[M^{'}]}\x\sigma )^{'}a,x_{0})\mid\cr\cr
&=&\mid\tau ^{'}\mid\mid(a,\sigma^{"}x_{0})\mid}$$
\p\vs 2mm
o\`u $\mid\tau ^{'}\mid$ d\'esigne le morphisme naturel :
\p\vs 2mm
$$\matrix{\mid\phi _{M^{'}}\mid&\Fhd{}{}&\mid\phi _{M}\mid\cr\cr
\mid(a,y)\mid&\Fhd{}{}&\mid\Bigl((\tau \circ I_{[m]})^{'}a,y\Bigr)\mid}$$
\p\vs 2mm
Donc l'image par $\phi$ devient :
\p\vs 2mm
$$\matrix{\phi \Bigl((\tau \x\sigma)^{'}a,(x_{0},x)\Bigl)&=&
\Bigl(\mid\tau ^{'}\mid\mid(a,\sigma^{"}x_{0})\mid,x\Bigl)\cr\cr
&\sim&\Bigl(\mid(a,\sigma ^{"}x_{0})\mid,\tau ^{"}x)\Bigl)\cr\cr
&=&\phi \Bigl(a,(\tau \x\sigma)^{"}(x_{0},x)\Bigl)}$$
\p\vs 2mm
On en d\'eduit donc que $\phi$ pr\'eserve les relation d'\'equivalences
et par suite induit
un morphisme $\mid \phi\mid\  :\  \mid\Phi \mid\ \fhd{}{}\
\mid{\cal R}\Phi \mid$.
Le morphisme $\mid \phi\mid$ admet un inverse $\mid G\mid$
d\'efini pour tout
$(a,x_{0},x)$ dans $\Phi _{M,m}\x R^{m}\x R^{M}$ par :
\p
$$\mid G\mid \mid\Bigl(\mid(a,x_{0})\mid,x\Bigr)\mid:=
\mid(a,(x_{0},x))\mid$$
\p
Ce qui montre que $\mid \phi\mid$ est un hom\'eomorphisme,
et par cons\'equent
les deux espaces ${\cal R}^{n}\Phi $ et $\mid\Phi \mid$ sont
hom\'eomorphes.
\p\hfill\fin
\p\vs 4mm\hs 4mm
{\bf Remarque (3.2.3):} Lorsqu'on rempla\c ce dans la proposition
pr\'ec\'edente $\Phi $ par
$\pi_{0}(\Omega^{n}X )$ on voit que ${\cal R}^{n}\pi_{0}(\Omega^{n}X )$
est hom\'eomorphe \`a $\mid\pi_{0}(\Omega^{n}X)\mid$, lequel est
hom\'eomorphe \`a $\mid\Pi_{n}(X)\mid$, car les deux $n$-groupo\"{i}des
$\pi_{0}(\Omega^{n}X)$ et $\Pi_{n}(X)$ sont isomorphes. Donc on obtient un
hom\'eomorphisme :
\p
$${\cal R}^{n}\pi_{0}(\Omega^{n}X)\ \Fhd{}{}\	\mid\Pi_{n}(X)\mid $$
\p\vs 4mm\hs 4mm
{\bf Proposition (3.2.4) :} {\it Soit $X$ un espace topologique. Alors
pour tout $n\geq 1$
il existe une \'equivalence d'homotopie faible :
\p
$${\cal R}^{n}\Omega ^{n}X\ \Fhd{\mu ^{n}_{X}}{}\ X$$
\p
compatible avec la fonctorialit\'e en $X$.}
\p\vs 4mm\hs 4mm
{\bf Preuve :} Le th\'eor\`eme (3.1.4) montre que la proposition est
vraie pour $n=1$.
Supposons maintenant qu'elle soit vraie jusqu'\`a un rang $n$.
Soit $m$ un objet de
$\Delta $, donc en appliquant l'hypoth\`ese de r\'ecurrence \`a l'espace
$Z(m)=(\Omega X)_{m}$, il existe une e.h.f :
\p
$${\cal R}^{n}\Omega ^{n}Z(m)\ \Fhd{\mu ^{n}_{Z(m)}}{}\ Z(m)$$
fonctorielle en $Z(m)$. Or on sait que $Z=\Omega X$ est un foncteur, alors :
\p
$${\cal R}^{n}\Omega ^{n+1}X\ \Fhd{\mu ^{n}_{Z}}{}\ \Omega X$$
\p
est un morphisme dans ETS. Donc d'apr\`es le Th\'eor\`eme (3.1.5)
on obtient une
e.h.f :
\p
$${\cal R}^{n+1}\Omega ^{n+1}X\ \Fhd{{\cal R}\mu ^{n}_{Z}}{}\
{\cal R}\Omega X$$
\p
Et en composant avec le morphisme $\mu ^{1}_{X} : {\cal R}\Omega X\fhd{}{}X$,
on
obtient notre morphisme $\mu ^{n+1}_{X}$ qui est fonctoriel en $X$ car les
foncteurs
$\Omega $ et ${\cal R}$ le sont.
\hfill\fin\p\vs 4mm
{\bf (3.3).--- LES ESPACES $\mid\Pi_{n}(X)\mid$ ET $X$ DANS $n$-Ho-{\cal T}
op :}
\p\vs 4mm\hs 4mm
Soit $X$ un espace topologique. La surjection canonique $p : X\fhd{}{}
\pi_{0}(X)$
n'est pas continue en g\'en\'eral, si on munit $\pi_{0}(X)$ de la topologie
discr\`ete; exemple : si $X=A\cup B$ avec
\p\vs 3mm
\centerline{$A=\{(0,y)\ \mid\ -1\leq y\leq 1\}$\hs 1cm et \hs 1cm
$B=\{(x,sin(1/x))\ \mid \ 0 < x\leq 1/2\pi\}$}
\p\vs 3mm
L'espace $X$, munit de la topologie induite de celle du plan r\'eel $R^{2}$,
est connexe
mais n'est pas connexe par arcs et on a $\pi_{0}(X)=\{A,B\}$. L'application
$p$ n'est pas continue car $p(X)=\pi_{0}(X)$ n'est pas connexe.
\p\hs 3mm
Comme on ne peut pas passer de $X$ \`a $\pi_{0}(X)$ par une application
continue,
il est possible de faire appelle \`a un espace interm\'ediaire qui aura
cette propri\'et\'e. Soit $X$ un espace topologique on lui associe
fonctoriellement un espace $F(X)$ d\'efini par le produit fibr\'e :
\p
$$F(X):= X{\x}_{\pi_{0}(X)^{g}}\pi_{0}(X)^{d}$$
\p
o\`u $(\ )^{g}$ et $(\ )^{d}$ d\'esignent respectivement la topologie
grossi\`ere et la topologie discr\`ete. On peut aussi voir cet espace
comme la r\'eunion disjointe :
\p
$$F(X):= \displaystyle\coprod _{\alpha\in K}^{disc}{\cal C}_{\alpha}$$
\p
o\`u ${\cal C}_{\alpha}$ est une composante connexe par arcs et $K$
un syst\`eme de repr\'esentants fix\'e.
\p
On r\'ecup\`ere alors les deux projections comme deux applications
continue
vers $X$ et $\pi_{0}(X)$ :
\p
$$X\ \Fhg{Pr_{1}}{}\ F(X)\ \Fhd{Pr_{2}}{}\ \pi_{0}(X)$$
\p\vs 4mm\hs 4mm
{\bf Proposition (3.3.1) :} {\it Pour tout espace topologique $X$,
le morphisme
$Pr_{1}$ est une \'equivalence d'homotopie faible. Si en plus $X$
est 0-tronqu\'e le morphisme $Pr_{2}$ est aussi une \'equivalence
d'homotopie faible.}
\p\vs 4mm\hs 4mm
{\bf Preuve :} Un \'el\'ement de $F(X)$ s'\'ecrit par $(x,\alpha)$
avec $\alpha$ dans
$K$ et $x$ dans la composante ${\cal C}_{\alpha}$. Si $\gamma$ est
un chemin dans
$F(X)$ entre $(x,\alpha)$ et $(y,\beta)$, alors $Pr_{1}\circ\gamma$
est un chemin
dans $X$ entre $x$ et $y$, ce qui montre que $\alpha =\beta$. On en
d\'eduit que
$Pr_{1}$ induit une bijection sur les $\pi_{0}$.
Soient $(x,\alpha)$ un \'el\'ement de $F(X)$, $i\geq 1$ un entier et
$\Gamma : (S^{i},s_{0})\fhd{}{}(F(X),(x,\alpha))$ un morphisme.
Le morphisme $Pr_{2}\circ\Gamma$ est \`a valeurs dans l'espace discret
$\pi_{0}(X)$,
donc constant (car $S^{i}$ est connexe) et prend la valeur
${\cal C}_{\alpha}$.
Ce qui montre que $\Gamma $ est compl\`etement d\'etermin\'e par
sa projection
$Pr_{1}\circ\Gamma$. Par suite $Pr_{1}$ induit des isomorphismes
sur les $\pi_{i}$.
Donc $Pr_{1}$ est une e.h.f.
Lorsque $X$ est 0-tronqu\'e, les groupes $\pi_{i}(X,x)$ sont nuls
pour $i\geq 1$,
il en est de m\^eme des groupes $\pi_{i}(F(X),(x,\alpha))$.
D'autre part
pour $i\geq 1$ les groupes $\pi_{i}(\pi_{0}(X),{\cal C}_{\alpha})$
sont nuls et
$\pi_{0}(\pi_{0}(X))$ est en bijection avec $\pi_{0}(X)$. D'o\`u
$\Pr_{2}$ est une e.h.f.
\hfill\fin\p\vs 4mm\hs 4mm
La correspondance $X\fhd{}{}F(X)$ est en fait fonctorielle et par
cons\'equent donne un foncteur
$F$ : $n$-ETS $\fhd{}{}$ $n$-ETS qui envoie $\Phi$ vers $F\circ\Phi$.
Par suite on obtient
deux morphismes :
\p
$$\Phi \ \Fhg{\alpha (\Phi ) }{}\ F\circ\Phi \ \Fhd{\beta (\Phi )}
{}\ \pi_{0}\circ\Phi $$
\p
avec pour tout $M$ objet de $\Delta ^{n}$, $\alpha (\Phi _{M})$ est
une \'equivalence
d'homotopie faible et $\beta (\Phi )_{M}$ en est une lorsque $\Phi _{M}$
est 0-tronqu\'e.
\p\vs 4mm\hs 4mm
{\bf Proposition (3.3.2) :} {\it Soit $\Phi$ un objet de $n$-ETS tel
que pour tout
objet $M$ de $\Delta$ l'espace $\Phi_{M}$ est 0-tronqu\'e. Alors il
existe deux
\'equivalences d'homotopies faibles, qui sont fonctorielles en $\Phi $ :
\p
$${\cal R}^{n}\Phi \ \Fhg{{\cal R}^{n}\alpha (\Phi ) }{}\ {\cal R}^{n}(F\circ
\Phi )\ \Fhd{{\cal R}^{n}\beta (\Phi )}{}\ {\cal R}^{n}(\pi_{0}\circ\Phi )$$}
\p\hs 4mm
{\bf Preuve :} Montrons la proposition par r\'ecurrence sur $n$.
\p
Pour $n=1$ il suffit d'appliquer le Th\'eor\`eme (3.1.5) aux deux
morphismes :
\p
$$\Phi \ \Fhg{\alpha (\Phi ) }{}\ F\circ\Phi \ \Fhd{\beta (\Phi )}
{}\ \pi_{0}\circ\Phi $$
\p
qui donne alors deux \'equivalences d'homotopies faibles :
\p
$${\cal R}\Phi \ \Fhg{{\cal R}\alpha(\Phi ) }{}\ {\cal R}(F\circ
\Phi )\ \Fhd{{\cal R}\beta (\Phi )}{}\ {\cal R}(\pi_{0}\circ\Phi)$$
\p
Supposons la proposition vraie jusqu'\`a un rang $n$, et soit $\Phi $ un
objet de ($n$+1)-ETS. En appliquant l'hypoth\`ese de r\'ecurrence
\`a $\Phi _{m}$
et en utilisant la fonctorialit\'e de $\Phi $ en $m$
on obtient deux morphismes dans la cat\'egorie 1-ETS :
\p
$${\cal R}^{n}\Phi \ \Fhg{{\cal R}^{n}\alpha (\Phi ) }{}\ {\cal R}^{n}
(F\circ\Phi)\ \Fhd{{\cal R}^{n}\beta (\Phi )}{}\ {\cal R}^{n}
(\pi_{0}\circ\Phi )$$
\p
et puis par le Th\'eor\`eme (3.1.5) on obtient :
\p
$${\cal R}^{n+1}\Phi \ \Fhg{{\cal R}^{n+1}\alpha (\Phi ) }{}
\ {\cal R}^{n+1}(F\circ
\Phi)\ \Fhd{{\cal R}^{n+1}\beta (\Phi )}{}\ {\cal R}^{n+1}
(\pi_{0}\circ\Phi )$$
\p\hfill\fin
\p\vs 4mm\hs 4mm
{\bf Proposition (3.3.3) :} {\it Soit $\Phi $ un objet de $n$-ETS tel
que pour tout objet $M$ de $\Delta ^{n}$ l'espace $\Phi _{M}$ soit
$i$-tronqu\'e avec $i\geq 1$. Alors pour tout $N$ objet de $\Delta ^{n+1}$
l'espace $\Omega\Phi _{_{N}}$ est ($i$-1)-tronqu\'e.}
\p\vs 4mm\hs 4mm
{\bf Preuve :} Soient $\Phi $ un objet de $n$-ETS et $(M,m)$ un objet de
$\Delta ^{n}\x\Delta $. L'ensemble $\Phi _{M}$ est un espace toplogique, donc
$\Omega\Phi _{_{M}}$ est un espace simplicial de Segal et en particulier on a
une e.h.f :
\p\vs 4mm
\centerline{$(\Omega\Phi _{_{M}})_{_{m}}\Fhd{\delta_{[m]}}{}
(\Omega \Phi _{_{M}})_{_{1}}{\x}_{(\Phi _{_{M}})^{disc}}
\cdots{\x}_{(\Phi _{_{M}})^{disc}}(\Omega \Phi _{_{M}})_{_{1}}$}
\p\vs 3mm
Donc il suffit de montrer que $\Psi _{1}:= (\Omega \Phi _{_{M}})_{_{1}}$
est ($i$+1)-tronqu\'e.
\p
 Remarquons d'abord que :
\p\vs 4mm
\centerline{$\Psi _{1}= \displaystyle\coprod _{x,y\in\Phi _{M}}Ch^{x,y}
(\Phi _{M})$
\hs 1cm et \hs 1cm $\Psi _{1}(x,y)=Ch^{x,y}(\Phi _{M})$}
\p\vs 4mm
La composition des chemins permet de construire, chaque fois qu'on a un
 chemin $g$
dans $\Psi _{1}(x,y)$, une \'equivalence d'homotopie :
\p	\vs 3mm
\centerline{$\Psi _{1}(x,y) \Fhd{}{}Lacet^{x}(\Phi _{M})$\hs 4mm o\`u
\hs 4mm
$Lacet^{x}(\Phi _{_{M}})=Ch^{x,x}(\Phi _{_{M}})$}
\p\vs 3mm
Soit $f$ un \'el\'ement de $\Psi _{1}(x,y)$, alors
$Lacet^{f}(\Psi _{1})=Lacet^{f}(\Psi _{1}(x,y)$, ce qui montre qu'on a
une
\'equivalence d'homotopie :
\p\vs 3mm
\centerline{$ Lacet^{f}(\Psi _{1}) \Fhd{}{} Lacet^{gf}\Bigl(Lacet^{x}
(\Phi _{M})\Bigr)$}
\p\vs 3mm
En particulier, pour tout entier $k\geq 1$ le groupe $\pi_{k}
(\Psi _{1},f)$ est
isomorphe \`a $\pi_{k+1}(\Phi _{M},x)$ et $\pi_{0}(\Psi _{1}(x,y),f)$
isomorphe \`a
$\pi_{1}(\Phi _{M},x)$.
\hfill\fin
\p\vs 4mm\hs 4mm
{\bf Remarque (3.3.4) :} Lorsqu'on pose $\Phi =\Omega^{n}X$ dans la
Proposition (3.3.2)
on obtient le diagramme d'\'equivalences d'homotopie faibles
suivant :
\p
\centerline{$\diagram{{\cal R}^{n}\Omega^{n}X&\Fhg{\sim}{}&
{\cal R}^{n}(F\circ\Omega^{n}X)&\Fhd{\sim}{}
&{\cal R}^{n}(\pi_{0}\circ\Omega^{n}X)
\cr \fvb{\wr}{}&&&&\fvb{}{\wr}\cr X&&&&\mid\Pi_{n}(X)\mid}$}
\p
ce qui montre qu'on a deux \'equivalences d'homotopie faibles,
fonctorielles en $X$ :
\p
\centerline{$\diagram{X&\Fhg{\sim}{}&{\cal F}_{n}(X)&
\Fhd{\sim}{\gamma }&\mid\Pi_{n}(X)\mid}$}
\p
o\`u  ${\cal F}_{n}(X)={\cal R}^{n}(F\circ\Omega^{n}X)$.
Les espaces ${\cal F}_{n}(X)$ et $\mid\Pi_{n}(X)\mid$ sont des
r\'ealisations
g\'eom\'etriques, donc des CW-complexes et par suite $\gamma $
est une \'equivalence
d'homotopie. Lorsqu'on regarde ce diagramme dans la cat\'egorie
$n$-Ho-{\cal T}op,
l'inverse de $\gamma $ permet d'obtenir un isomorphisme de
$\mid\Pi_{n}(X)\mid$ vers $X$ compatible avec la fonctorialit\'e
en $X$
\p\vs 4mm
{\bf (3.4).---  LES $n$-GROUPOIDES $\Phi $ ET $\Pi_{n}(\mid\Phi \mid)$
DANS $n$-Ho-Gr :}
\p\vs 4mm\hs 3mm
Dans ce paragraphe on va traiter la deuxi\`eme \'etape qui consiste \`a
construire pour tout $n$-groupo\"{i}de $\Phi $ une \'equivalence
ext\'erieure entre
$\Phi $ et $\Pi_{n}(\mid\Phi \mid)$ et qui va correspondre \`a
un isomorphisme
dans la cat\'egorie $n$-Ho-Gr.
Soit $\Phi $ un $n$-groupo\"{i}de. On peut construire de fa\c con
naturelle un morphisme
\p
$${\cal L} : \Phi \ \Fhd{}{}\ \Pi_{n}(\mid\Phi \mid)$$
\p
d\'efini pour tout objet $M$ de $\Delta ^{n}$ par :
\p
\centerline{$\diagram{\Phi _{M}&\Fhd{{\cal L}_{_{M}}}{}&
\overline{\mid\Phi \mid}_{M}
\cr a&\Fhd{}{}&{\cal L}_{_{M}}(a):=\overline{f}}$ \hs 1cm
avec \hs 1cm
$\diagram{R^{M}&\Fhd{f}{}&\mid\Phi \mid\cr
x&\Fhd{}{}&\mid(a,x)\mid}$}
\p
Montrons d'abord que ce morphisme est bien d\'efini c.\`a.d que
$f$ est un
\'el\'ement de ${\mid\Phi \mid}_{M}$. Soient  $(m_{1},\cdots,m_{n})$
un objet de
$\Delta ^{n}$ et $1\leq k < n$ un entier, et posons :
\p\vs 3mm
\centerline{$\matrix{N=(m_{1},\cdots,m_{k-1}) && S=(m_{k+1},\cdots,m_{n})
\cr\cr
\gamma  =I_{N}\x I_{[0]}\x\delta _{0}&&\gamma _{0}=I_{N}\x I_{[0]}\x
e _{0_{S}}\cr\cr
\alpha   =I_{N}\x \delta _{i}\x I_{S}&avec&0\leq i \leq m_{k}}$}
\p\vs 3mm
O\`u $\delta _{i}$ est l'application qui envoie 0 vers $i$ et
$e _{0_{S}}$ est l'unique application
de $S$ vers $0_{S}$. Soient les diagrammes :
\p
$$N\x [m_{k}]\x S\ \Fhg{\alpha }{}\ N\x [0]\x S\
\Fhd{\Fhg{\gamma }{}}{\gamma _{0}}\ N\x [0]\x 0_{S}$$
\p\vs 3mm
$$R^{S}\x R^{m_{k}}\x R^{N}\ \Fhg{\alpha ^{"}}{}\ R^{S}\x R^{0}\x R^{N}\
\Fhd{\Fhg{\gamma ^{"}}{}}{\gamma _{0}^{"}}\ R^{0_S}\x R^{0}\x R^{N}$$
\p\vs 3mm
$$\Phi _{N,m_{k},S}\ \Fhg{\alpha ^{'}}{}\ \Phi _{N,0,S}\
\Fhd{\Fhg{\gamma ^{'}}{}}{\gamma _{0}^{'}}\ \Phi _{N,0,0_{S}}$$
\p\vs 3mm
Comme $\Phi $ est un $n$-groupo\"{i}de, le foncteur $\Phi _{N}$ est constant
donc
$\Phi _{N,0,S}=\Phi _{N,0,0_{S}}$ et $\gamma ^{'}=\gamma _{0}^{'}=
id_{\Phi _{N,0,0_{S}}}$.
\p\vs 3mm
Soit $\Bigl( a, (x,\delta ^{"}_{i}(1),y)\Bigr)\in R^{S}\x R^{m_{k}}
\x R^{N}$ donc
$f(x,\delta ^{"}_{i}(1),y)=\mid\Bigl( a, (x,\delta ^{"}_{i}(1),y)
\Bigr)\mid$.
\p\vs 3mm
\centerline{$\matrix{\Bigl( a, (x,\delta ^{"}_{i}(1),y)\Bigr)&=&
\Bigl( a, \alpha ^{"}(x,1,y)\Bigr)\cr\cr&\sim&\Bigl( \alpha ^{'}a,
(x,1,y)\Bigr)\cr\cr
&=&\Bigl(\gamma ^{'}_{0}\gamma ^{'} \alpha ^{'}a, (x,1,y)\Bigr)\cr\cr
&\sim&\Bigl( \gamma ^{'} \alpha ^{'}a, \gamma ^{"}_{0}(x,1,y)\Bigr)\cr	\cr
&=&\Bigl( \gamma ^{'} \alpha ^{'}a, (1_{S},1,y)\Bigr)\cr\cr
&\sim&\Bigl( \alpha ^{'}a, \gamma ^{"}(1_{S},1,y)\Bigr)\cr\cr
&=&\Bigl(\alpha ^{'}a, (\delta ^{"}_{0_{S}}(1),1,y)\Bigr)\cr	\cr
&\sim&\Bigl(a, \alpha ^{"}(\delta ^{"}_{0_{S}}(1),1,y)\Bigr)\cr\cr
&=&\Bigl(a, (\delta ^{"}_{0_{S}}(1),\delta ^{"}_{i}(1),y)\Bigr)}$}
\p\vs 4mm
On en d\'eduit que $f(x,\delta ^{"}_{i}(1),y)=f(\delta ^{"}_{0_{S}}(1),
\delta ^{"}_{i}(1),y)$
ne d\'epend pas de $x$ et ce ci pour tout $1\leq k\leq n$, ce qui montre
que $f$ est
bien dans ${\mid\Phi \mid}_{M}$.
On veut prouver par r\'ecurrence sur $n$ que le morphisme ${\cal L}$
est une
$n$-\'equivalence ext\'erieure. Donc on aimerait bien savoir si une
$n$-\'equivalence ext\'erieure peut se d\'eduire d'une ($n$-1)-\'equivalence
ext\'erieure.
D'abord pour un $n$-nerf $\Phi $ et pour tout couple d'objets $(x,y)$
de $\Phi $,
on d\'esigne par $\Phi _{1}(x,y)$ la ($n$-1)-pr\'e-cat\'egorie
d\'efinie pour tout
objet $M$ de $\Delta ^{n-1}$ par :
\p
$$\Phi _{1}(x,y)(M):=\Phi _{1,M}(x,y)=\{\sigma \in \Phi _{1,M} \
\mid\ \delta ^{'}_{0}
(\sigma )=x\ et\ \delta ^{'}_{1}(\sigma )=y\}$$
\p
\centerline{avec \hs 1cm $\delta ^{'}_{i}=(\delta _{i}\x I_{_{M}})^{'} :
\Phi _{1,M}\ \Fhd{}{}\
\Phi _{0,M}=\Phi _{0_{n}}$}
\p\vs 4mm
Il est ais\'e de voir que $\Phi _{1}(x,y)$ est muni d'une structure de
($n$-1)-nerf
induite par celle de $\Phi _{1}$; c'est {\it le ($n$-1)-nerf des fl\`eches
de $\Phi $ qui ont
pour source et but respectivement $x$ et $y$}. Lorsque $\Phi _{1}$ est un
($n$-1)-groupo\"{i}de il en est de m\^eme pour $\Phi _{1}(x,y)$.
\p\vs 4mm\hs 4mm
{\bf Proposition (3.4.1) :} {\it Un morphisme $F : \Phi \fhd{}{}\Psi $ entre
$n$-groupo\"{i}des est une $n$-\'equivalence ext\'erieure si et seulement si :
\p
(i) l'application $\pi_{0}(F) : \pi_{0}(\Phi )\fhd{}{}\pi_{0}(\Psi )$ est
surjective.
\p
(ii) Pour tout $x,y$ dans Ob$(\Phi )$ :
$$F_{1}(x,y) : \Phi _{1}(x,y)\ \Fhd{}{}\ \Psi _{1}(F(x),F(y))$$
\p
est une ($n$-1)-\'equivalence ext\'erieure.}
\p\vs 4mm\hs 4mm
{\bf Preuve :} Remarquons d'abord que (ii) montre que $\pi_{0}(F)$ est aussi
injective.
\p
D'apr\`es la proposition (2.2.4)dans [11], une $n$-\'equivalence ext\'erieure
entre
deux groupo\"{i}des se traduit par des isomorphismes entre leurs groupes
d'homotopies.
\p
Pour tout $x,y$ objets de $\Phi $ et tout $f$ objet de ${\cal C}_{i}(\Phi )
(x,y)$ avec
$2\leq i\leq n$ on a :
\p\vs 3mm
\centerline{$\pi_{0}(\Phi _{1}(x,y))=T^{n-1}(\Phi _{1}(x,y))=(T^{n-1}
\Phi _{1})(x,y)=
Hom_{{\cal C}_{1}(\Phi )}(x,y)\simeq \pi_{1}(\Phi ,x)$}
\p\vs 3mm
\centerline{${\cal C}_{i-1}(\Phi _{1}(x,y))=T^{n-i}\Phi _{I_{(i-1)}}(x,y)=
{\cal C}_{i}(\Phi )(x,y)$}
\p\vs 3mm
\centerline{$\pi_{i}(\Phi ,f)=\pi_{i-1}(\Phi _{1}(x,y),f)$}
\p\vs 3mm
L'axiome (ii) est \'equivalent de dire que $F$ induit des isomorphismes sur
les $\pi_{i}$ pour
$1\leq i\leq n$, ce qui montre la proposition.
\hfill\fin
\p\vs 2mm\hs 4mm
Soit $\Phi $ une $n$-pr\'e-cat\'egorie. On d\'esigne par ${\cal P}\Phi $
l'espace topologique simplicial d\'efini pour tout $m$ par :
$({\cal P}\Phi )_{m}={\cal R}^{n-1}\Phi _{m}$. On a alors ${\cal R}
{\cal P}\Phi =
{\cal R}^{n}\Phi \simeq\mid\Phi \mid$.
\p\vs 4mm\hs 4mm
{\bf Proposition (3.4.2) :} {\it
(i) Soit $\Phi $ un $n$-groupo\"{i}de, alors ${\cal P}\Phi $ est un espace
simplicial de Segal
avec :
\p
\centerline{$T(\pi_{0}({\cal P}\Phi ))=T^{n}\Phi $\hs 4mm et \hs 4mm
$\pi_{0}({\cal P}\Phi )\simeq T^{n-1}\Phi $}
\p\vs 3mm
(ii) Soit $F:\Phi \fhd{}{}\Psi $ une $n$-\'equivalence ext\'erieure entre
$n$-groupo\"{i}des,
alors le morphisme ${\cal R}^{n}F : {\cal R}^{n}\Phi \fhd{}{}
{\cal R}^{n}\Psi $ est une
\'equivalence d'homotopie faible.}
\p\vs 4mm\hs 4mm
{\bf Preuve :} On va montrer la proposition par r\'ecurrence sur $n$.
\p
Cas $n=1$ : Un groupo\"{i}de $\Phi $ est en particulier un espace simplicial
de Segal
lorsqu'on consid\`ere la topologie discr\`ete sur les ensemble $\Phi _{m}$.
On a
${\cal P}\Phi =\Phi $ et $\pi_{0}(\Phi )=\Phi $ , d'o\`u
$T(\pi_{0}(\Phi ))=T\Phi $.
\p
Soit $F : \Phi \fhd{}{}\Psi $ une 1-\'equivalence
ext\'erieure entre groupo\"{i}des. Donc pour tous $x,y$ objets de $\Phi $
 on a :
\p
$$F_{1}(x,y) : \Phi _{1}(x,y)\ \Fhd{}{}\ \Psi  _{1}(x^{'},y^{'})$$
\p
est une bijection, consid\'er\'ee comme hom\'eomorphisme par la topologie
discr\`ete.
\p
D'apr\`es le Th\'eor\`eme (3.2.1) de Segal, on a deux e.h.f :
\p\vs 3mm
\centerline{$\Phi _{1}(x,y)\ \Fhd{}{}\ Ch^{x,y}({\cal R}\Phi )$}
\p\vs 3mm
\centerline{$\Psi  _{1}(x^{'},y^{'})\ \Fhd{}{}\ Ch^{x^{'},y^{'}}
({\cal R}\Psi )$}
\p\vs 3mm
qui rendent commutatif le diagramme suivant :
\p
\centerline{$\diagram{\Phi _{1}(x,y)&\Fhd{F_{1}(x,y)}{}&\Psi  _{1}
(x^{'},y^{'})\cr
\fvb{}{}&&\fvb{}{}\cr
Ch^{x,y}({\cal R}\Phi )&\Fhd{}{Ch^{x,y}({\cal R}F)}&Ch^{x^{'},y^{'}}
({\cal R}\Psi )}$}
\p
et qui induisent un diagramme commutatif sur les $\pi_{i}$ pour
$ i\geq 0$. Ce qui montre que
${\cal R}F$ induit des isomorphismes sur les $\pi_{i} $ pour
$i\geq 1$. D'autre part on a le
diagramme commutatif suivant :
\p
\centerline{$\diagram{T\Phi &\Fhd{\sim}{}&T\Psi  \cr
\fvb{\wr}{}&&\fvb{}{\wr}\cr
\pi_{0}({\cal R}\Phi )&\Fhd{}{\pi_{0}({\cal R}F)}&\pi_{0}({\cal R}\Psi )}$}
\p
qui entra\^{\i}ne que $\pi_{0}({\cal R}F)$ est une bijection et par
cons\'equent
${\cal R}F : {\cal R}\Phi \fhd{}{}{\cal R}\Psi $ est une e.h.f.
Supposons la proposition vraie pour $n$ et montrons qu'elle l'est
pour $n$+1.
\p
Soit $F : \Phi \fhd{}{}\Psi $ une ($n$+1)-\'equivalence ext\'erieure
entre
($n$+1)-groupo\"{i}des. Pour tout entier $m\geq 2$, on applique
l'hypoth\`ese de
r\'ecurrence \`a la $n$-\'equivalence :
\p
$$\Phi _{m}\ \Fhd{\delta _{[m]}}{}\ \Phi _{1}{\x}_{_{\Phi _{0}}}
\cdots{\x}_{_{\Phi _{0}}}\Phi _{1}$$
\p
on obtient l'e.h.f :
\p
$${\cal R}^{n}\Phi _{m}\ \Fhd{}{}\ {\cal R}^{n}(\Phi _{1}{\x}_{_{\Phi _{0}}}
\cdots{\x}_{_{\Phi _{0}}}\Phi _{1})$$
\p\vs 4mm\hs 4mm
{\bf Lemme (3.4.3) :} {\it Si $\Phi $ est un espace topologique
$n$-simplicial,
alors pour tout entier $m\geq 2$ le morphisme $\theta (\Phi )$
dans le
diagramme commutatif suivant :
\p
\centerline{$\diagram{{\cal R}^{n-1}\Phi _{m}&\Fhd{{\cal R}^{n-1}
\delta _{[m]}}{}&
{\cal R}^{n-1}(\Phi _{1}{\x}_{_{\Phi _{0}}}\cdots{\x}_{_{\Phi _{0}}}
\Phi _{1})
\cr \Vert&&	\fvb{}{\theta (\Phi )}\cr {\cal R}^{n-1}\Phi _{m}&
\Fhd{}{\delta _{[m]}}&
{\cal R}^{n-1}\Phi _{1}{\x}_{_{{\cal R}^{n-1}\Phi _{0}}}
\cdots{\x}_{_{{\cal R}^{n-1}
\Phi _{0}}}{\cal R}^{n-1}\Phi _{1}}$}
est une \'equivalence d'homotopie faible (e.h.f).}
\p\vs 4mm\hs 4mm
{\bf Preuve :} On sait que si $K$ et $L$ sont deux ensemble simpliciaux,
l'application continue $\theta  : {\cal R}(K\x L)\fhd{}{}{\cal R}(K)
\x{\cal R}(L)$ poss\`ede
une application inverse qui est continue sur toute cellule produit de
${\cal R}(K)\x{\cal R}(L)$
( J. Peter May [2]). Cette situation montre que $\theta $ induit des
isomorphismes sur
les groupes d'homotopies et on conclut que $\theta $ est une
e.h.f. Un raisonnement analogue, qui tient compte de la topologie,
permet de prolonger
cela \`a un produit fibr\'e fini d'espaces topologiques simpliciaux.
On en d\'eduit que le Lemme est vrai pour $n=2$. Supposons le Lemme
vrai
pour $n\geq 2$ et montrons qu'il est vrai pour $n$+1.
Soit $\Phi $ un espace topologique ($n$+1)-simplicial simplicial.
En appliquant l'hypoth\`ese de r\'ecurrence \`a l'espace topologique
simplicial
${\cal R}\Phi$, on obtient une e.h.f :
\p
$${\cal R}^{n-1}({\cal R}\Phi _{1}{\x}_{_{{\cal R}\Phi _{0}}}
\cdots{\x}_{_{{\cal R}
\Phi _{0}}}{\cal R}\Phi _{1})\ \Fhd{\theta ({\cal R}\Phi)}{}\
{\cal R}^{n}\Phi _{1}{\x}_{_{{\cal R}^{n}\Phi _{0}}}
\cdots{\x}_{_{{\cal R}^{n}\Phi _{0}}}{\cal R}^{n}\Phi _{1}$$
\p
D'autre part, pour tout $M$ objet de $\Delta ^{n-1}$, on applique
l'hypoth\`ese de
r\'ecurrence \`a l'ensemble simplicial $\Phi _{1,M}$ on obtient une e.h.f :
\p
$${\cal R}(\Phi _{1,M}{\x}_{_{\Phi _{0_{n-1}}}}\cdots{\x}_{_{\Phi _{0_{n-1}}}}
\Phi _{1,M})\ \Fhd{\mu_{_{M}}}{}\ {\cal R}\Phi _{1,M}{\x}_{_{{\cal R}
\Phi _{0_{n}}}}
\cdots{\x}_{_{{\cal R}\Phi _{0_{n}}}}{\cal R}\Phi _{1,M}$$
\p
Comme $\mu_{_{M}}$ est fonctoriel en $M$, on lui applique le
Th\'eor\`eme
(3.1.5)
(Bousfield-Kan [BK]) ($n$-1)-fois, ce qui donne une e.h.f :
\p
$${\cal R}^{n}(\Phi _{1}{\x}_{_{\Phi _{0}}}\cdots{\x}_{_{\Phi _{0}}}
\Phi _{1})\ \Fhd{{\cal R}^{n}(\mu)}{}\ {\cal R}^{n-1}
({\cal R}\Phi _{1}{\x}_{_{{\cal R}\Phi _{0}}}
\cdots{\x}_{_{{\cal R}\Phi _{0}}}{\cal R}\Phi _{1})$$
\p\vs 2mm
Le compos\'e $\theta ({\cal R}\Phi)\circ{\cal R}^{n}(\mu)=\theta (\Phi )$, et
par cons\'equent $\theta (\Phi )$ est  une e.h.f
\hfill\fin\p\vs 2mm\hs 4mm
D'apr\`es le Lemme (3.4.3) et en remarquant que ${\cal R}^{n}
\Phi_{m} =({\cal P}\Phi)_{m}$
on voit que :
\p
$$({\cal P}\Phi)_{m}\ \Fhd{}{}\
({\cal P}\Phi) _{1}{\x}_{_{({\cal P}\Phi )_{0}}}
\cdots{\x}_{_{({\cal R}\Phi )_{0}}}({\cal P}\Phi )_{1}$$
\p
est une e.h.f. On en d\'eduit que ${\cal P}\Phi $
v\'erifie l'axiome (ii) de Segal.
\p\vs 2mm
L'espace $({\cal P}\Phi )_{0}\simeq\mid\Phi _{0}\mid=\Phi _{0_{n+1}}$ car
le foncteur
$\Phi _{0}$ est constant, d'o\`u l'axiome (i) de Segal.
\p
Montrons maintenant que :
$T(\pi_{0}({\cal P}\Phi ))=T^{n+1}\Phi $ , et pour cela on pose :
\p\vs 3mm
\centerline{${\cal A}=\pi_{0}({\cal P}\Phi )$\hs 1cm et \hs 1cm
${\cal B}=\pi_{0}({\cal P}(T\Phi ))$}
\p\vs 3mm
D'apr\`es l'hypoth\`ese de r\'ecurrence ${\cal P}(T\Phi )$ est un espace
simplicial
de Segal, donc ${\cal B}$ est un groupo\"{i}de et $T{\cal B}=
T^{n+1}\Phi $.
\p\vs 3mm
\centerline{$\matrix{{\cal R}^{n}\Phi _{0}=
[\displaystyle\coprod _{m}{\cal R}^{n-1}\Phi _{0,m}\x R^{m}]^{\sim}=
{\cal R}^{n-1}\Phi _{0,0}\hs 4mm car \hs 4mm(\Phi _{0,m}=\Phi _{0,0})\cr\cr
alors \hs 1cm T{\cal A}=[{\cal R}^{n}\Phi _{0}]^{\sim}=[{\cal R}^{n-1}
\Phi _{0,0}]^{\sim}=
T{\cal B}=T^{n+1}\Phi }$}
\p \vs 3mm
En vertu du Th\'eor\`eme de Segal on a les isomorphismes :
\p\vs 3mm
\centerline{$\matrix{{\cal A}(m)=\pi_{0}({\cal R}{\cal P}
\Phi _{m})&\Fhd{\sim}{}&
T[\pi_{0}({\cal P}\Phi _{m})]}$}
\p\vs 3mm
Et par hypoth\`ese de r\'ecurrence on a :
\p\vs 3mm
\centerline{$T[\pi_{0}({\cal P}\Phi _{m})]\ =\ T^{n}\Phi _{m}$}
\p\vs 3mm
On en d\'eduit que ${\cal A}\simeq T^{n}\Phi $ et par suite
${\cal P}\Phi $
satisfait les axiomes (i) et (iii) de Segal (c'est alors un
espace simplicial de Segal).
Montrons maintenant que le morphisme
\p
${\cal R}^{n+1}F : {\cal R}^{n+1}\Phi \fhd{}{}{\cal R}^{n+1}\Psi $
est une
e.h.f.
\p
Comme $F : \Phi \fhd{}{}\Psi $ est une ($n$+1)-\'equivalence ext\'erieure,
cela implique
d'apr\`es la Proposition (4.3.1) que pour tout $x,y$ dans Ob$(\Phi )$ :
\p
$$F_{1}(x,y) : \Phi _{1}(x,y)\ \Fhd{}{}\ \Psi _{1}(x^{'},y^{'})$$
\p
est une $n$-\'equivalence ext\'erieure, et par l'hypoth\`ese de
r\'ecurrence une e.h.f :
\p
$${\cal R}^{n}F_{1}(x,y) : {\cal R}^{n}\Phi _{1}(x,y)\ \Fhd{}{}\
{\cal R}^{n}\Psi _{1}(x^{'},y^{'})$$
\p
D'apr\`es le Th\'eor\`eme de Segal et la Proposition (3.4.2),
appliqu\'es \`a ${\cal P}\Phi $, on a d'une part un diagramme
commutatif :
\p
\centerline{$\diagram{T^{n}\Phi &=&T(\pi_{0}({\cal P}\Phi ))&
\Fhd{\sim}{}&
\pi_{0}({\cal R}^{n+1}\Phi )\cr
\fvb{\wr}{}&&&&\fvb{}{\pi_{0}({\cal R}^{n+1}F)}\cr
T^{n}\Psi&=&T(\pi_{0}({\cal P}\Psi ))&\Fhd{\sim}{}&
\pi_{0}({\cal R}^{n+1}\Psi )}$}
\p
qui entra\^{\i}ne que $\pi_{0}({\cal R}^{n+1}F)$ est une
bijection.
\p
Et d'autre part, comme ${\cal R}^{n}\Phi _{1}(x,y)=({\cal P}\Phi)_{1}(x,y)$,
le morphisme $F :\Phi \fhd{}{}\Psi $ induit un morphisme
${\cal P}\Phi  : {\cal P}\Phi \fhd{}{}{\cal P}\Psi $ et puis un diagramme
commutatif :
\p
\centerline{$\diagram{({\cal P}\Phi)_{1}(x,y)&\Fhd{e.h.f}{}&
({\cal P}\Psi )_{1}(x^{'},y^{'})
\cr\fvb{e.h.f}{}&&\fvb{}{e.h.f}\cr
Ch^{x,y}({\cal R}^{n+1}\Phi )&\Fhd{}{}&Ch^{x^{'},y^{'}}({\cal R}^{n+1}\Psi )}$}
\p
lequel induit des isomorphismes sur les $\pi_{i}$ pour $i\geq 0$, ce
qui montre que
${\cal R}^{n+1}F$ induit des isomorphismes sur les $\pi_{i}$ pour $i\geq 1$.
Par cons\'equent ${\cal R}^{n+1}F$ est une e.h.f.
\p\hfill\fin
\p\vs 4mm\hs 4mm
{\bf Proposition (3.4.4) :} {\it Soient $f :X\fhd{}{}Y$ une \'equivalence
d'homotopie faible,
$n\geq 1$ un entier et $\Phi $ un $n$-groupo\"{i}de. Alors on a :
\p
(i) Le morphisme $\Pi_{n}(f): \Pi_{n}(X)\ \fhd{}{}\ \Pi_{n}(Y)$ est une
$n$-\'equivalence
ext\'erieure.
\p
(ii) Le morphisme ${\cal L} : \Phi \ \fhd{}{}\	\Pi_{n}(\mid \Phi \mid)$ est
une $n$-\'equivalence ext\'erieure.}
\p\vs 4mm\hs 4mm
{\bf Preuve :} On va effectuer pour cela une r\'ecurrence sur $n$.
\p
\item{\bf (a)}{\bf Cas n=1} : (i) Si $f : X\fhd{}{}Y$ est une e.h.f alors elle
induit une bijection sur
entre les $\pi_{0}$ et des isomorphismes entre les $\pi_{i}$ pour $i\geq 1$.
D'apr\`es le Th\'eor\`eme (2.4.4) du chapitre 2 , pour tout $x$ dans $X$
il existe un isomorphisme de groupes :
\p
$$\pi_{1}(\Pi_{1}(X),x)\ \Fhd{\sim}{}\	\pi_{1}(X,x)$$
\p
fonctoriel en $X$. Et comme $T(\Pi_{1}(X))=\pi_{0}(X)$, on d\'eduit que
$\Pi_{1}(f)$ est une
1-e.h.f.
\p
(ii) Un groupo\"{i}de $\Phi $ est en particulier un espace simplicial de
Segal munit de la
topologie discr\`ete. Donc d'apr\`es le Th\'eor\`eme de Segal on a :
\p
$$T\Phi =T(\pi_{0}\circ\Phi )\ \Fhd{\sim}{}\ \pi_{0}(\mid\Phi\mid)=
T(\Pi_{1}(\mid\Phi \mid))$$
\p
est une bijection et pour tout objets $x,y$ de $\Phi $ on a une e.h.f :
\p
$$\Phi _{1}(x,y)\ \Fhd{}{}\ Ch^{x,y}(\mid\Phi \mid)$$
\p
et puis une bijection
\p
$$\Phi _{1}(x,y)=\pi_{0}(\Phi _{1}(x,y))\ \Fhd{\sim}{}
\ \pi_{0}(Ch^{x,y}(\mid\Phi \mid))=\Pi_{1}(\mid\Phi \mid)_{_{1}}(x,y)$$
\p
ce qui montre d'apr\`es la Proposition (4.3.1) que
${\cal L} : \Phi \ \fhd{}{}\	\Pi_{1}(\mid \Phi \mid)$ est une
$1$-\'equivalence
ext\'erieure.
\p
\item{\bf (b)} Supposons la Proposition vraie pour $n$
et montrons qu'elle l'est pour $n$+1.
\p
(ii) Soit $\Phi $ un ($n$+1)-groupo\"{i}de. On sait que ${\cal P}\Phi $ est
un espace de Segal
donc pour tout couple $(x,y)$ d'objets de $\Phi $ on a un isomorphisme et
une e.h.f :
\p
$$\mid\Phi _{1}(x,y)\mid\ \fhd{\sim}{}\ ({\cal P}\Phi)_{1}(x,y)
\ \Fhd{e.h.f}{}(\Omega\mid\Phi \mid )_{1}(x,y)$$
\p
Soit $\lambda $ le morphisme compos\'e des deux morphismes ci-dessus,
il est d\'efini comme suite :
\p\vs 3mm
\centerline{$\forall (a_{_{M}},u_{_{M}}) \in \Phi _{1,M}(x,y)\x R^{M}$\hs 3mm
avec \hs 3mm
$M\in Ob(\Delta ^{n})$ \hs 3mm on pose \hs 3mm $w=\mid(a_{_{M}},u_{_{M}})
\mid$}
\p\vs 5mm
\centerline{donc \hs 1cm $\matrix{\lambda (w) : R\ 	\Fhd{}{}\ \mid \Phi
\mid (x,y)
\cr\cr t	\ \Fhd{}{}\ \mid(a_{_{M}},(u_{_{M}},t))\mid}$}
\p\vs 3mm
Consid\'erons maintenant le diagramme dans la cat\'egorie $n$-Gr suivant :
\p
\centerline{$\diagram{ \Phi _{1}(x,y)&\Fhd{{\cal L}^{'}}{}&
\Pi_{n}(\mid\Phi _{1}(x,y)\mid)\cr\fvb{{\cal L}_{1}(x,y)}{}&(D)&\fvb{}
{\Pi_{n}(\lambda )}
\cr \Pi_{n}(\mid\Phi \mid)_{_{1}}(x,y)&\Fhg{\sim}{F_{1}}&
\Pi_{n}((\Omega\mid\Phi \mid)_{_{1}})(x,y)}$}
\p
Dans le diagramme ci-dessus les morphismes ${\cal L}^{'}$, $\Pi_{n}(\lambda )$
sont
des $n$-\'equivalences ext\'erieures par hypoth\`ese de r\'ecurrence et $F_{1}$
est un
isomorphisme d'apr\`es la preuve de la proposition (3.2.1). En explicitant ces
trois morphismes on montre que le diagramme (D) est commutatif, ce qui
entra\^ine que
${\cal L}_{1}(x,y)$ est une $n$-\'equivalence ext\'erieure.
\p
Le diagramme $(D)$ est commutatif si et seulement si pour tout objet $M$
de $\Delta ^{n}$
le diagramme $(D)_{_{M}}$ suivant est commutatif :
\p
\centerline{$\diagram{ \Phi _{1,M}(x,y)&\Fhd{{\cal L}^{'}_{M}}{}&
\overline{\mid\Phi _{1}(x,y)\mid}_{M}\cr
\fvb{{\cal L}_{1,M}(x,y)}{}&(D)_{_{M}}&\fvb{}{\Pi_{n}(\lambda )_{M}}
\cr \overline{\mid\Phi \mid}_{_{1,M}}(x,y)&\Fhg{}{F_{1,M}}&
\overline{(\Omega\mid\Phi \mid)}_{_{1,M}}(x,y)}$}
\p
\centerline{Soit \hs 1cm$\tau \in \Phi _{1,M}(x,y)$ \hs 1cm donc \hs 1cm
$F_{1,M}\circ\Pi_{n}(\lambda )_{M}\circ{\cal L}^{'}_{M}(\tau )=
\overline{\Gamma }$}
\p\vs 5mm
\centerline{o\`u \hs 1cm $\matrix{\Gamma  : &R^{M}\x R&\Fhd{}{}&
\mid\Phi \mid\cr\cr
&(a,t)&\Fhd{}{}&\mid(\tau ,(a,t))\mid}$}
\p\vs 5mm
Donc $\overline{\Gamma }={\cal L}_{1,M}(x,y)(\tau )$. Ce qui montre
que $(D)_{_{M}}$ est
commutatif et par suite $(D)$ est commutatif.
\p
Par ailleurs on a un diagramme commutatif :
\p
\centerline{$\diagram{T^{n+1}\Phi &\Fhd{T^{n+1}{\cal L}}{}&\pi_{0}
(\mid\Phi \mid)\cr
\Vert&&\fvb{}{\wr}\cr T(\pi_{0}({\cal P}\Phi ))&\Fhd{}{\sim}&\pi_{0}
({\cal R}({\cal P}\Phi ))}$}
car ${\cal P}\Phi $ est un espace de Segal et ${\cal R}({\cal P}\Phi) $
hom\'eomorphe \`a
$\mid\Phi \mid$. Donc ${\cal L}$ est une ($n$+1)-\'equivalence
ext\'erieure.
\p
(i) Soit $f :X\fhd{}{}Y$ une e.h.f. Donc pour tout $x,y$ dans $X$
on a une bijection $\pi_{0}(f) : \pi_{0}(X)\fhd{}{}\pi_{0}(Y)$ et
une e.h.f :
\p\vs 4mm
\centerline{$(\Omega X)_{_{1}}(x,y)\ \Fhd{h}{}\ (\Omega Y)_{_{1}}
(x^{'},y^{'})$
\hs 4mm o\`u \hs 4mm $x^{'}=f(x)$ et $y^{'}=f(y)$}
\p\vs 3mm
Par application de l'hypoth\`ese de r\'ecurrence on obtient une
$n$-\'equivalence ext\'erieure :
\p
$$\Pi_{n}((\Omega X)_{_{1}}(x,y))\ \Fhd{}{}\ \Pi_{n}((\Omega Y)_{_{1}}
(x^{'},y^{'}))$$
\p
Consid\'erons le diagramme suivant  :
\p
\centerline{$\diagram{\Pi_{n}((\Omega X)_{_{1}}(x,y))&\Fhd{}{}&
\Pi_{n}((\Omega Y)_{_{1}}(x^{'},y^{'}))\cr\fvb{\wr}{}&(E)&\fvb{}{\wr}\cr
\Pi_{n+1}(X)_{_{1}}(x,y)&\Fhd{}{\Pi_{n+1}({\cal L})_{_{1}}(x,y)}&
\Pi_{n+1}(Y)_{_{1}}(x^{'},y^{'})}$}
\p
o\`u les fl\`eches verticales sont des isomorphismes. La commutativit\'e
de ce diagramme permet de conclure que $\Pi_{n+1}({\cal L})_{_{1}}(x,y)$
est une $n$-\'equivalence ext\'erieure.
\p
Soit $M$ un objet de $\Delta ^{n}$, on va di\'ecrire le diagramme $(E)_{_{M}}$
suivant par deux chemins diff\'erents :
\p
\centerline{$\diagram{\overline{(\Omega X)}_{_{1,M}}(x,y)&\Fhd{}{}&
\overline{(\Omega Y)}_{_{1,M}}(x^{'},y^{'})\cr\fvb{\wr}{}&(E)_{_{M}}&
\fvb{}{\wr}\cr
\overline{X}_{_{1,M}}(x,y)&\Fhd{}{}&\overline{Y}_{_{1,M}}(x^{'},y^{'})}$}
\p
\p
\centerline{$\diagram{\overline{(\Omega X)}_{_{1,M}}(x,y)&\Fhd{}{}&
\overline{(\Omega Y)}_{_{1,M}}(x^{'},y^{'})&\Fhd{}{}&\overline{Y}_{_{1,M}}
(x^{'},y^{'})\cr
\overline{\sigma }&\Fhd{}{}&\overline{\sigma ^{'}}&\Fhd{}{}&
\overline{\sigma ^{"}}}$}
\p
Les repr\'esentants $\sigma ^{'}$ et $\sigma ^{"}$ sont d\'efinis
comme suit :
\p\vs 3mm
\centerline{$\sigma ^{'}=h\circ\sigma $ \hs  4mm et \hs 4mm
$\sigma ^{"}(a,t) =\sigma ^{'}(a)(t)$\hs 4mm
pour tout  $(a,t)$  dans  $R^{M}\x R$}
\p\vs 3mm
par l'autre chemin on a :
\p
\centerline{$\diagram{\overline{(\Omega X)}_{_{1,M}}(x,y)&\Fhd{}{}&
\overline{X}_{_{1,M}}(x,y)&\Fhd{}{}&\overline{Y}_{_{1,M}}(x^{'},y^{'})\cr
\overline{\sigma }&\Fhd{}{}&\overline{\sigma _{1}}&\Fhd{}{}&
\overline{\sigma _{2}}}$}
\p
o\`u cette fois les repr\'esentants $\sigma _{1}$ et $\sigma _{2}$ sont
d\'efinis par :
\p\vs 3mm
\centerline{$\sigma _{2}=f\circ\sigma _{1}$ \hs  4mm et \hs 4mm
$\sigma _{1}(a,t) =\sigma (a)(t)$\hs 4mm
pour tout  $(a,t)$  dans  $R^{M}\x R$}
\p\vs 3mm
\centerline{Or 	\hs 4mm
$\sigma ^{"}(a,t)=\sigma ^{'}(a)(t)=h[\sigma (a)](t)=f[\sigma (a)(t)]=
f[\sigma _{1}(a)(t)]=\sigma _{2}(a,t)$}
\p\vs 3mm
Donc le diagramme $(E)_{_{M}}$ est commutatif et puis $(E)$ est commutatif.
D'autre part la bijection $\pi_{0}(f)$  co\"{i}ncide avec l'application :
\p
$$ T^{n+1}[\Pi_{n+1}(f)]\  :\ T^{n+1}[\Pi_{n+1}(X)]\ \Fhd{}{}
\ T^{n+1}[\Pi_{n+1}(Y)]$$
\p
On en d\'eduit alors d'apr\`es la Proposition (3.4.1) que $\Pi_{n+1}(f)$
est une
($n$+1)-\'equivalence ext\'erieure.
\p\hfill\fin\p\vs 4mm
{\bf Conclusion et preuve du Th\'eor\`eme (3.0):}
\p\vs 4mm\hs 3mm
Au chapitre 2 nous avons construit les deux foncteurs :
\p
\centerline{$\diagram{n-Top\ \Fhd{\Pi_{n}(\ )}{\Fhg{}{\mid\ \mid}}\ n-Gr}$}
\p
r\'ealisation g\'eom\'etrique $\mid\ \mid$ et $n$-groupo\"{i}de de
Poincar\'e $\Pi_{n}(\ )$. Au cours de ce chapitre nous avons pu
\'etudier leur comportement et on peut conclure les r\'esultats suivant :
\p
(i) Si $f : X\fhd{}{}Y$ est une e.h.f alors
$\Pi_{n}(f) : \Pi_{n}(X)\fhd{}{}\Pi_{n}(Y)$ est une $n$-\'equivalence
ext\'erieure.
\p
(ii) Si $\phi : \Phi\fhd{}{}\Phi$ est une $n$-est une $n$-\'equivalence, alors
$\mid\phi\mid : \mid\Phi\mid\fhd{}{}\mid\Phi\mid$ est une e.h.f
\p
ce qui permet de dire que les foncteurs ci-dessus induisent des foncteurs sur
les cat\'egories localis\'ees :
\p\vs 2mm
\centerline{$\diagram{n-Ho-Top\ \Fhd{\Pi_{n}(\ )}{\Fhg{}{\mid\ \mid}}
\ n-Ho-Gr}$}
\p\vs 2mm
D'autre part on a des isomorphismes naturels :
\p\vs 2mm
\centerline{$\matrix{\Phi&\Fhd{}{}&\Pi_{n}(\mid\Phi\mid)}$ \hs 4mm dans
 $n$-Ho-Gr}
\p\vs 2mm
\centerline{$\matrix{\mid\Pi_{n}(X)\mid&\Fhd{}{}&X}$ \hs 4mm dans
$n$-Ho-top}
\p\vs 2mm
qui sont compatibles avec la fonctorialit\'e en $\Phi$ et en $X$
respectivement.
Donc les foncteurs $\Pi_{n}(\ )$ et $\mid\ \mid$ induisent une \'equivalence
entre les cat\'egories
$n$-Ho-Gr et $n$-Ho-top.
\p\hfill\fin\p\vs 1cm
{\bf  Remerciements:}
\p\vs 4mm\hs 4mm
{\it Je tiens \`a exprimer ma profonde gratitude au Professeur Carlos SIMPSON
qui a dirig\'e ce travail avec grand soin, alternant une aide efficace \`a de
pr\'ecieux conseils. J'exprime aussi mes sinc\`eres remerciements aux
Professeurs J. TAPIA et J-C. SIKORAV pour leur conseils et
l'int\'er\^et qu'ils portent \`a ce travail.}
\p\vfill\eject
\centerline{\bf REFERENCES BIBLIOGRAPHIQUES}
\vs 1cm
\item{\bf [As]} {\bf  J. F Adams :} {\it Infinite loop spaces,}
{ Herman weyl lectures,
Annals of mathematics studies N 90. Princeton university press. 1978.}
\p\vs 2mm
\item{\bf [BK]} {\bf  A. K. Bousfield-D. M. Kan :} {\it Homotopy limits
completions and localisations,}
{ LNM 304, Springer-Verlag  1972.}\p\vs 2mm
\item{\bf [ES]} {\bf  Edwin-H.Spanier :} {\it Algebraic Topology,}
{ Mc Graw-Hill Book Company, 1966.}\p\vs 2mm
\item{\bf [GZ]} {\bf  P. Gabriel-M. Zisman :} {\it Calculus of
fractions and homotopy theory,}
{ Ergebnisse der mathematik, Vol. 35. Berlin-Heidelberg-New York :
Springer-Verlag (1967).}\p\vs 2mm
\item{\bf [JR]} {\bf Joseh J. Rotman :} {\it An introduction to
algebraic topologgy,}
{ Graduate texts in mathematics, 119 New York Berlin Heidelberg
London Paris,
Springer-Verlag, 1988.}\p\vs 2mm
\item{\bf [M]} {\bf S. Mac Lane :} {\it Categories for the working
mathematicians,} {Spinger-Verlag,
New York (1971).}
\item{\bf [MM]} {\bf  S. Mac Lane-I. Moerdijk :} {\it Sheaves in
geometry and logic,} { Springer-Verlag,
New York (1992).}\p\vs 2mm
\item{\bf [My]} {\bf  J. P. May :} {\it Simplicial objects in
algebraic topology,}{ Princeton : Van Nostrand
mathematical studies 11, (1967).}\p\vs 2mm
\item{\bf [RS]} {\bf  Robert M. Switzer :} {\it  Alrebraic
topology-Homotopy and Homology,}
{ Die Grundlehren der mathematischen wissenschaften
Springer-Verlag Berlin Heidelberg New York, 1975.}\p\vs 2mm
\item{\bf [SH]} {\bf  Sze-Tsen Hu :} {\it Homotopy theory,}
{ Wayne state university,
Detroit Michigan, Academic press.  New York and London, 1959.}\p\vs 2mm
\item{\bf [Sl]} {\bf  Graeme Segal :} {\it Homotopy everything H-spaces,}
 {AMS.}\p\vs 2mm
\item{\bf [Ti]} {\bf  Z. Tamsamani :} {\it On non strict notions of
$n$-caterory and
$n$-groupo\"{i}de via multi-simplicial sets,} {AMS Preprint Server
(a.k.a. alg-geom/9512006).}\p\vs 2mm
\end